\documentclass{amsart}
\usepackage{amsmath,amssymb}
\usepackage{yhmath}

\newcommand{\bdis}{\begin{displaymath}}
\newcommand{\edis}{\end{displaymath}}
\newcommand{\be}{\begin{equation}}
\newcommand{\ee}{\end{equation}}
\newcommand{\mbb}{\mathbb}
\newcommand{\mcal}{\mathcal} 
\newcommand{\mfrak}{\mathfrak}

\newcommand{\vp}{\varphi}

\newcommand{\zf}{\zeta\left(\frac{1}{2}+it\right)}

\newcommand{\FR}{\frac{x^n+y^n}{z^n}}

\theoremstyle{definition}

\theoremstyle{remark}
\newtheorem{remark}[]{Remark}

\newtheorem*{mydef11}{{\bf Theorem 1}}

\newtheorem*{mydef12}{{\bf Theorem 2}}

\newtheorem*{mydef13}{{\bf Theorem 3}}

\newtheorem*{mydef51}{{\bf Lemma 1}}

\newtheorem*{mydef52}{{\bf Lemma 2}}

\newtheorem*{mydef53}{{\bf Lemma 3}}

\newtheorem*{mydef54}{{\bf Lemma 4}}

\newtheorem*{mydef81}{{\bf Property 1}}

\newtheorem*{mydef82}{{\bf Property 2}}

\newtheorem*{mydefII1}{\bf Consequence 1}

\numberwithin{equation}{section}

\begin{document}

\title[Jacob's ladders, next equivalents \dots]{Jacob's ladders, next equivalents of the Fermat-Wiles theorem and new infinite sets of the equivalents generated by the Dirichlet's series}

\author{Jan Moser}

\address{Department of Mathematical Analysis and Numerical Mathematics, Comenius University, Mlynska Dolina M105, 842 48 Bratislava, SLOVAKIA}

\email{jan.mozer@fmph.uniba.sk}

\keywords{Riemann zeta-function}

\begin{abstract}
In this paper we obtain new sets of equivalents of the Fermat-Wiles theorem. Simultaneously, we obtain also asymptotic connections between the set of Dirichlet's series, certain segments of the Dirichlet's sum $\mfrak{D}(x)$, Riemann zeta-function and Jacob's ladders. 
\end{abstract}
\maketitle

\section{Introduction} 

In this paper, that continues our series \cite{7} -- \cite{15}, we obtain new points of contact between the Riemann's zeta-function and the Fermat-Wiles theorem. 

\subsection{} 

First, we use simultaneously: 
\begin{itemize}
	\item[(a)] Hardy-Littlewood formula (1922) 
	\be \label{1.1} 
	\int_1^T|\zeta(\sigma+it)|{\rm d}t\sim\zeta(2\sigma)T,\ T\to\infty, 
	\ee 
	\item[(b)] Selberg's formula 
	\be \label{1.2} 
	\int_1^T|S_1(t)|^{2l}{\rm d}t\sim \bar{c}(L)T, 
	\ee  
	where 
	\bdis 
	S_1(t)=\frac{1}{\pi}\int_0^t\arg\left\{\zeta\left(\frac 12+iu\right)\right\}{\rm d}u, 
	\edis  
	\item[(c)] and our almost linear formula 
	\be \label{1.3} 
	\int_T^{\overset{1}{T}}\left|\zf\right|^2{\rm d}t\sim(1-c)T;\ \overset{1}{T}=\vp_1^{-1}(T), 
	\ee  
	for every fixed 
	\bdis 
	\epsilon>0,\ \sigma\geq\frac 12+\epsilon,\ l\in\mbb{N}.
	\edis 
\end{itemize}

On this basis, we obtain, for example, the following result: The $\zeta$-condition 
\be \label{1.4} 
\begin{split}
& \lim_{\tau\to\infty}\frac{1}{\tau}\left\{
\int_{\FR\frac{\tau}{\sqrt[3]{(1-c)\bar{c}(l)\zeta(2\sigma)}}}^{[\FR\frac{\tau}{\sqrt[3]{(1-c)\bar{c}(l)\zeta(2\sigma)}}]^1}\left|\zf\right|^2{\rm d}t\times \right. \\
& \left. \int_1^{\FR\frac{\tau}{\sqrt[3]{(1-c)\bar{c}(l)\zeta(2\sigma)}}}|\zeta(\sigma+it)|^2{\rm d}t \times \right. \\ 
& \left. \int_1^{\FR\frac{\tau}{\sqrt[3]{(1-c)\bar{c}(l)\zeta(2\sigma)}}}|S_1(t)|^{2l}{\rm d}t
\right\}\not=1 
\end{split}
\ee 
on the set of all Fermat's rationals represents the next $\zeta$-equivalent of the Fermat-Wiles theorem for every fixed 
\bdis 
\epsilon>0,\ \sigma\geq\frac 12+\epsilon,\ l\in\mbb{N}, 
\edis 
where, of course, 
\bdis 
[Y]^1=\vp_1^{-1}(Y),\ x,y,z,n\in\mbb{N},\ n\geq 3. 
\edis 

\begin{remark}
We have obtained also a similar result for certain linear combination of the basic almost linear integrals (\ref{1.1}) -- (\ref{1.3}). 
\end{remark} 

\subsection{} 

Now, let us remind the old mean-value theorem for Dirichlet's series: if 
\be \label{1.5} 
f(s)=\sum_{n=1}^\infty \frac{a_n}{n^s},\ s=\sigma+it 
\ee 
is absolutely convergent for $\sigma=\sigma_0$, then for $\sigma=\sigma_0$, 
\be \label{1.6} 
\lim_{T\to\infty}\frac{1}{T}\int_0^T|f(\sigma_0+it)|^2{\rm d}t=\sum_{n=1}^\infty\frac{|a_n|^2}{n^{2\sigma_0}}. 
\ee  

\begin{remark}
Of course, the Riemann's zeta-function 
\bdis 
\zeta(s)=\sum_{n=1}^\infty \frac{1}{n^s},\ \sigma>1 
\edis 
represents the canonical element of the infinite set of all Dirichlet's series $\mfrak{D}$. 
\end{remark} 

Let 
\be \label{1.7} 
\sum_{n=1}^\infty\frac{|a_n|^2}{n^{2\sigma_0}}=F(\sigma_0;f), 
\ee  
where 
\be \label{1.8} 
F(\sigma_0;f)=F[\sigma_0(f);f]\in\mbb{R}^+
\ee 
for every fixed element $f\in\mfrak{D}$ and corresponding point of absolute convergence $\sigma_0(t)$. Then we have\footnote{Comp. \cite{14}, (1.28).} the following functional 
\be \label{1.9} 
\lim_{\tau\to\infty}\frac{1}{\tau}\int_0^{\frac{x}{F(\sigma_0;f)}\tau}|f(\sigma_0+it)|^2{\rm d}t=x,\ x>0. 
\ee  

\subsection{} 

In this paper we shall consider the finite sets 
\be \label{1.10} 
f_m(\sigma_0^m+it)=\sum_{n=1}^\infty\frac{a_n^m}{n^{\sigma_0^m+it}},\ m=1,\dots,M,\ M\in\mbb{N}
\ee 
of Dirichlet's series, where 
\be \label{1.11} 
\sigma_0^m=\sigma_0^m(f_m),\ m=1,\dots,M 
\ee 
for every fixed 
\bdis 
M\in\mbb{N},\ \{f_m\}_{m=1}^M,\ \{\sigma_0^m\}_{m=1}^M. 
\edis 
In this direction, we obtain the following result: the $\mfrak{D}$-condition 
\be \label{1.12} 
\lim_{\tau\to\infty}\frac{1}{\tau}\prod_{m=1}^M\left\{\int_0^{\FR\frac{\tau}{\{\prod_m F_m(\sigma_0^m)\}^{1/M}}}\right\}^{\frac{1}{M}}|f_m(\sigma_0^m+it)|^2{\rm d}t\not=1
\ee 
on the set of all Fermat's rationals expresses the next $\mfrak{D}$-equivalent of the Fermat-Wiles theorem for every fixed 
\bdis 
M\in\mbb{N},\ \{f_m\}_{m=1}^M,\ \{\sigma_0^m\}_{m=1}^M. 
\edis 

\begin{remark}
Of course, the formula (\ref{1.12}) represents the infinite set of of $\mfrak{D}$-equivalents. 
\end{remark}

\section{Jacob's ladders: notions and basic geometrical properties}  

\subsection{} 

In this paper we use the following notions of our works \cite{2} -- \cite{6}: 
\begin{itemize}
\item[{\tt (a)}] Jacob's ladder $\vp_1(T)$, 
\item[{\tt (b)}] direct iterations of Jacob's ladders 
\bdis 
\begin{split}
	& \vp_1^0(t)=t,\ \vp_1^1(t)=\vp_1(t),\ \vp_1^2(t)=\vp_1(\vp_1(t)),\dots , \\ 
	& \vp_1^k(t)=\vp_1(\vp_1^{k-1}(t))
\end{split}
\edis 
for every fixed natural number $k$, 
\item[{\tt (c)}] reverse iterations of Jacob's ladders 
\be \label{2.1}  
\begin{split}
	& \vp_1^{-1}(T)=\overset{1}{T},\ \vp_1^{-2}(T)=\vp_1^{-1}(\overset{1}{T})=\overset{2}{T},\dots, \\ 
	& \vp_1^{-r}(T)=\vp_1^{-1}(\overset{r-1}{T})=\overset{r}{T},\ r=1,\dots,k, 
\end{split} 
\ee   
where, for example, 
\be \label{2.2} 
\vp_1(\overset{r}{T})=\overset{r-1}{T}
\ee  
for every fixed $k\in\mbb{N}$ and every sufficiently big $T>0$. We also use the properties of the reverse iterations listed below.  
\be \label{2.3}
\overset{r}{T}-\overset{r-1}{T}\sim(1-c)\pi(\overset{r}{T});\ \pi(\overset{r}{T})\sim\frac{\overset{r}{T}}{\ln \overset{r}{T}},\ r=1,\dots,k,\ T\to\infty,  
\ee 
\be \label{2.4} 
\overset{0}{T}=T<\overset{1}{T}(T)<\overset{2}{T}(T)<\dots<\overset{k}{T}(T), 
\ee 
and 
\be \label{2.5} 
T\sim \overset{1}{T}\sim \overset{2}{T}\sim \dots\sim \overset{k}{T},\ T\to\infty.   
\ee  
\end{itemize} 

\begin{remark}
	The asymptotic behaviour of the points 
	\bdis 
	\{T,\overset{1}{T},\dots,\overset{k}{T}\}
	\edis  
	is as follows: at $T\to\infty$ these points recede unboundedly each from other and all together are receding to infinity. Hence, the set of these points behaves at $T\to\infty$ as one-dimensional Friedmann-Hubble expanding Universe. 
\end{remark}  

\subsection{} 

Let us remind that we have proved\footnote{See \cite{8}, (3.4).} the existence of almost linear increments 
\be \label{2.6} 
\begin{split}
& \int_{\overset{r-1}{T}}^{\overset{r}{T}}\left|\zf\right|^2{\rm d}t\sim (1-c)\overset{r-1}{T}, \\ 
& r=1,\dots,k,\ T\to\infty,\ \overset{r}{T}=\overset{r}{T}(T)=\vp_1^{-r}(T)
\end{split} 
\ee 
for the Hardy-Littlewood integral (1918) 
\be \label{2.7} 
J(T)=\int_0^T\left|\zf\right|^2{\rm d}t. 
\ee  

For completeness, we give here some basic geometrical properties related to Jacob's ladders. These are generated by the sequence 
\be \label{2.8} 
T\to \left\{\overset{r}{T}(T)\right\}_{r=1}^k
\ee 
of reverse iterations of of the Jacob's ladders for every sufficiently big $T>0$ and every fixed $k\in\mbb{N}$. 

\begin{mydef81}
The sequence (\ref{2.8}) defines a partition of the segment $[T,\overset{k}{T}]$ as follows 
\be \label{2.9} 
|[T,\overset{k}{T}]|=\sum_{r=1}^k|[\overset{r-1}{T},\overset{r}{T}]|
\ee 
on the asymptotically equidistant parts 
\be \label{2.10} 
\begin{split}
& \overset{r}{T}-\overset{r-1}{T}\sim \overset{r+1}{T}-\overset{r}{T}, \\ 
& r=1,\dots,k-1,\ T\to\infty. 
\end{split}
\ee 
\end{mydef81} 

\begin{mydef82}
Simultaneously with the Property 1, the sequence (\ref{2.8}) defines the partition of the integral 
\be \label{2.11} 
\int_T^{\overset{k}{T}}\left|\zf\right|^2{\rm d}t
\ee 
into the parts 
\be \label{2.12} 
\int_T^{\overset{k}{T}}\left|\zf\right|^2{\rm d}t=\sum_{r=1}^k\int_{\overset{r-1}{T}}^{\overset{r}{T}}\left|\zf\right|^2{\rm d}t, 
\ee 
that are asymptotically equal 
\be \label{2.13} 
\int_{\overset{r-1}{T}}^{\overset{r}{T}}\left|\zf\right|^2{\rm d}t\sim \int_{\overset{r}{T}}^{\overset{r+1}{T}}\left|\zf\right|^2{\rm d}t,\ T\to\infty. 
\ee 
\end{mydef82} 

It is clear, that (\ref{2.10}) follows from (\ref{2.3}) and (\ref{2.5}) since 
\be \label{2.14} 
\overset{r}{T}-\overset{r-1}{T}\sim (1-c)\frac{\overset{r}{T}}{\ln \overset{r}{T}}\sim (1-c)\frac{T}{\ln T},\ r=1,\dots,k, 
\ee  
while our eq. (\ref{2.13}) follows from (\ref{2.6}) and (\ref{2.5}). 

\section{The equivalent of the Fermat-Wiles theorem generated by means of product of three simplest $\zeta$-integrals} 

\subsection{} 

We use the following simplest\footnote{Almost linear in the variable $T$.} $\zeta$-formulas: 
\begin{itemize}
	\item[(A)] Hardy-Littlewood formula (1922) 
	\be \label{3.1} 
	\int_1^T|\zeta(\sigma+it)|^2{\rm d}t=\zeta(2\sigma)T+\mcal{O}(T^{1-\epsilon}\ln T),\ T\to\infty 
	\ee  
	for every fixed\footnote{See \cite{17}, pp. 30, 31.} 
	\be \label{3.2} 
	\epsilon>0,\ \sigma\geq \frac{1}{2}+\epsilon, 
	\ee 
	\item[(B)] Selberg's formula 
	\be \label{3.3} 
	\int_1^T|S_1(t)|^{2l}{\rm d}t=\bar{c}(l)T+\mcal{O}\left(\frac{T}{\ln T}\right) 
	\ee  
	for every fixed $l\in\mbb{N}$ (see \cite{16}, p. 130), where 
	\be \label{3.4} 
	S_1(t)=\frac{1}{\pi}\int_1^t\arg\left\{\zf\right\}{\rm d}t, 
	\ee 
	\item[(C)] and also our formula 
	\be \label{3.5} 
	\begin{split} 
	& \int_{\overset{r-1}{T}}^{\overset{r}{T}}\left|\zf\right|^2{\rm d}t=(1-c)\overset{r-1}{T}+\mcal{O}(T^{1/3+\delta}), \\ 
	& r=1,\dots,k,\ T\to\infty 
	\end{split} 
	\ee 
	in the case $r=1$, for every fixed $k\in\mbb{N}$ and every fixed small $\delta>0$, where 
	\be \label{3.6} 
	\overset{r}{T}=\overset{r}{T}(T)=\vp_1^{-r}(T). 
	\ee 
\end{itemize} 

\begin{remark}
Hardy and Littlewood have proved the nontrivial part of the formula (\ref{3.1}), namely for the part of the critical strip and also for the line $\sigma=1$. 
\end{remark}

Consequently, we use the asymptotic forms of the formulas (\ref{3.1}), (\ref{3.3}) and (\ref{3.5}): 
\be \label{3.7} 
\begin{split}
& \int_1^T|\zeta(\sigma+it)|^2{\rm d}t\sim \zeta(2\sigma)T,\ \sigma\geq\frac 12+\epsilon, \\ 
& \int_1^T|S_1(t)|^{2l}{\rm d}t\sim \bar{c}(l)T, \\ 
& \int_{T}^{\overset{1}{T}}\left|\zf\right|^2{\rm d}t\sim (1-c)T, 
\end{split}
\ee 
for $T\to\infty$. 

\subsection{} 

Next, the basic formula for this section is the following one: 
\be \label{3.8} 
\begin{split}
& \left\{
\int_T^{\overset{1}{T}}\left|\zf\right|^2{\rm d}t \times \int_1^T|\zeta(\sigma+it)|^2{\rm d}t \times \int_1^T|S_1(t)|^{2l}{\rm d}t
\right\}^{1/3}\sim \\ 
& \sqrt[3]{(1-c)\bar{c}(l)\zeta(2\sigma)}T,\ T\to\infty. 
\end{split}
\ee 
Now, we use the following substitution 
\be \label{3.9} 
T=\frac{x}{\sqrt[3]{(1-c)\bar{c}(l)\zeta(2\sigma)}}\tau,\ \{T\to\infty\}\ \Leftrightarrow \ \{\tau\to\infty\}
\ee 
for every fixed $x>0$, into the eq. (\ref{3.8}) and we obtain the following functional. 

\begin{mydef51}
\be \label{3.10} 
\begin{split}
& \lim_{\tau\to\infty}\frac{1}{\tau}\left\{
\int_{\frac{x\tau}{\sqrt[3]{(1-c)\bar{c}(l)\zeta(2\sigma)}}}^{[\frac{x\tau}{\sqrt[3]{(1-c)\bar{c}(l)\zeta(2\sigma)}}]^1}\left|\zf\right|^2{\rm d}t\times \right. \\
& \left. \int_1^{\frac{x\tau}{\sqrt[3]{(1-c)\bar{c}(l)\zeta(2\sigma)}}}|\zeta(\sigma+it)|^2{\rm d}t \times \right. \\ 
& \left. \int_1^{\frac{x\tau}{\sqrt[3]{(1-c)\bar{c}(l)\zeta(2\sigma)}}}|S_1(t)|^{2l}{\rm d}t
\right\}^{1/3}=x 
\end{split}
\ee 
for every fixed 
\bdis 
x>0,\ l\in\mbb{N},\ \sigma\geq\frac 12+\epsilon. 
\edis  
\end{mydef51} 

Next, in the special case 
\bdis 
x\to\FR, 
\edis  
it follows from (\ref{3.10}) that the following holds true. 

\begin{mydef52}
\be \label{3.11} 
\begin{split} 
& \lim_{\tau\to\infty}\frac{1}{\tau}\left\{
\int_{\FR\frac{\tau}{\sqrt[3]{(1-c)\bar{c}(l)\zeta(2\sigma)}}}^{[\FR\frac{\tau}{\sqrt[3]{(1-c)\bar{c}(l)\zeta(2\sigma)}}]^1}\left|\zf\right|^2{\rm d}t\times \right. \\
& \left. \int_1^{\FR\frac{\tau}{\sqrt[3]{(1-c)\bar{c}(l)\zeta(2\sigma)}}}|\zeta(\sigma+it)|^2{\rm d}t \times \right. \\ 
& \left. \int_1^{\FR\frac{\tau}{\sqrt[3]{(1-c)\bar{c}(l)\zeta(2\sigma)}}}|S_1(t)|^{2l}{\rm d}t
\right\}^{1/3}=\FR 
\end{split} 
\ee  
for every fixed Fermat's rational and every fixed 
\bdis 
l\in\mbb{N},\ \sigma\geq \frac 12+\epsilon. 
\edis  
\end{mydef52}

Consequently, we have the following result.  

\begin{mydef11}
The $\zeta$-condition 
\be \label{3.12} 
\begin{split}
& \lim_{\tau\to\infty}\frac{1}{\tau}\left\{
\int_{\FR\frac{\tau}{\sqrt[3]{(1-c)\bar{c}(l)\zeta(2\sigma)}}}^{[\FR\frac{\tau}{\sqrt[3]{(1-c)\bar{c}(l)\zeta(2\sigma)}}]^1}\left|\zf\right|^2{\rm d}t\times \right. \\
& \left. \int_1^{\FR\frac{\tau}{\sqrt[3]{(1-c)\bar{c}(l)\zeta(2\sigma)}}}|\zeta(\sigma+it)|^2{\rm d}t \times \right. \\ 
& \left. \int_1^{\FR\frac{\tau}{\sqrt[3]{(1-c)\bar{c}(l)\zeta(2\sigma)}}}|S_1(t)|^{2l}{\rm d}t
\right\}^{1/3}\not=1 
\end{split}
\ee 
on the set of all Fermat's rationals represents the next $\zeta$-equivalent of the Fermat-Wiles theorem for every fixed 
\bdis 
l\in\mbb{N},\ \sigma\geq \frac 12+\epsilon. 
\edis 
\end{mydef11} 

\section{An equivalent defined by means of certain linear combination the three almost linear $\zeta$-integrals} 

\subsection{} 

The basic formula for this section is the following one\footnote{See (\ref{3.7}).} 
\be \label{4.1} 
\begin{split}
& \frac{1}{1-c}\int_T^{\overset{1}{T}}\left|\zf\right|^2{\rm d}t+\frac{1}{\zeta(2\sigma)}\int_1^T|\zeta(\sigma+it)|^2{\rm d}t+ \\ 
& \frac{1}{\bar{c}(l)}\int_1^T|S_1(t)|^{2l}{\rm d}t\sim 3T,\ T\to\infty 
\end{split}
\ee 
for every fixed 
\bdis 
l\in\mbb{N},\ \sigma\geq \frac 12+\epsilon. 
\edis 

Now, the substitution 
\be \label{4.2} 
T=\frac{x}{3}\tau 
\ee  
into (\ref{4.1}), for every fixed $x>0$, gives us the following functional. 

\begin{mydef53}
\be \label{4.3} 
\begin{split}
& \lim_{\tau\to\infty}\frac{1}{\tau}\left\{
\frac{1}{1-c}\int_{\frac{x}{3}\tau}^{[\frac{x}{3}\tau]^1}\left|\zf\right|^2{\rm d}t+\frac{1}{\zeta(2\sigma)}\int_1^{\frac{x}{3}\tau}|\zeta(\sigma+it)|^2{\rm d}t+ \right. \\ 
& \left. \frac{1}{\bar{c}(l)}\int_1^{\frac{x}{3}\tau}|S_1(t)|^{2l}{\rm d}t
\right\}=x
\end{split}
\ee 
for every fixed 
\be \label{4.4} 
x>0,\ l\in\mbb{N},\ \sigma\geq\frac{1}{2}+\epsilon. 
\ee 
\end{mydef53} 

Of course, in the special case 
\bdis 
x\to\FR
\edis 
we obtain from (\ref{4.3}) the following result. 

\begin{mydef12}
The $\zeta$-condition 
\be \label{4.5} 
\begin{split}
	& \lim_{\tau\to\infty}\frac{1}{\tau}\left\{
	\frac{1}{1-c}\int_{\FR\frac{\tau}{3}}^{[\FR\frac{\tau}{3}]^1}\left|\zf\right|^2{\rm d}t+\frac{1}{\zeta(2\sigma)}\int_1^{\FR\frac{\tau}{3}}|\zeta(\sigma+it)|^2{\rm d}t+ \right. \\ 
	& \left. \frac{1}{\bar{c}(l)}\int_1^{\FR\frac{\tau}{3}}|S_1(t)|^{2l}{\rm d}t
	\right\}\not=1
\end{split}
\ee 
on the class of all Fermat's rationals represents the next $\zeta$-equivalent of the Fermat-Wiles theorem for every fixed 
\bdis 
l\in\mbb{N},\ \sigma\geq\frac 12+\epsilon. 
\edis 
\end{mydef12} 

\section{New infinite set of equivalents of the Fermat-Wiles theorem generated by the Dirichlet's series without Jacob's ladders}

\subsection{} 

Let us remind that for any Dirichel's series there is the following corresponding formula\footnote{Comp. (\ref{1.6}) and (\ref{1.7}).}
\be \label{5.1} 
\int_0^T|f(\sigma_0+it)|^2{\rm d}t=F(\sigma_0;f)T+o(T),\ T\to\infty 
\ee 
for every fixed 
\be \label{5.2} 
f\in\mfrak{D},\ \sigma_0=\sigma_0(f), 
\ee 
Consequently, to the finite set of Dirichlet's series (\ref{1.10}) there is a set of corresponding formulas 
\be \label{5.3} 
\int_0^T|f_m(\sigma_0^m+it)|^2{\rm d}t\sim F_m(\sigma_0^m)T,\ T\to\infty, 
\ee 
where\footnote{Comp. (\ref{1.7}).} 
\be \label{5.4} 
F_m(\sigma_0^m)=F_m(\sigma_0^m;f_m)=\sum_{n=1}^\infty \frac{|a_n|^2}{n^{2\sigma_0^m}}\in\mbb{R}^+. 
\ee 
Now, we have 
\be \label{5.5} 
\prod_{m=1}^M\int_0^T|f_m(\sigma_0^m+it)|^2{\rm d}t\sim \left\{\prod_{m=1}^MF_m(\sigma_0^m)\right\}T^M, 
\ee 
i. e. 
\be \label{5.6} 
\begin{split}
& \prod_{m=1}^M\left\{\int_0^T|f_m(\sigma_0^m+it)|^2{\rm d}t\right\}^{1/M}\sim \left\{\prod_{m=1}^MF_m(\sigma_0^m)\right\}^{1/M}T, \\ 
& T\to\infty ;\ \prod_{m=1}^MF_m(\sigma_0^m)=\prod_{m=1}^MF_m(\sigma_0^m;f_m). 
\end{split}
\ee 
By making use of the substitution 
\be \label{5.7} 
T=\frac{x}{\left\{\prod_{m=1}^MF_m(\sigma_0^m)\right\}^{1/M}}\tau,\ x>0
\ee 
we obtain the following functional. 

\begin{mydef54}
\be \label{5.8} 
\begin{split}
& \lim_{\tau\to\infty}\frac{1}{\tau}\prod_{m=1}^M\left\{
\int_0^{\frac{x}{\left\{\prod_{m=1}^MF_m(\sigma_0^m)\right\}^{1/M}}\tau}|f_m(\sigma_0^m+it)|^2{\rm d}t
\right\}^{1/M}=x
\end{split}
\ee 
for every fixed 
\bdis 
M\in\mbb{N},\ f_m\in\mfrak{D},\ \sigma_0^m(f_m). 
\edis 
\end{mydef54} 

And finally, our last result implies the following equivalent holds true. 

\begin{mydef13}
The $\mfrak{D}$-condition 
\be \label{5.9} 
\lim_{\tau\to\infty}\frac{1}{\tau}\prod_{m=1}^M\left\{
\int_0^{\FR\frac{\tau}{\left\{\prod_{m=1}^MF_m(\sigma_0^m)\right\}^{1/M}}}|f_m(\sigma_0^m+it)|^2{\rm d}t
\right\}^{1/M}\not=1
\ee 
on the set of all Fermat's rationals represents the next $\mfrak{D}$-equivalent of the Fermat-Wiles theorem for every fixed 
\bdis 
M\in\mbb{N},\ f_m(\sigma_0^m+it)\in\mfrak{D},\ \sigma_0^m(f_m),\ m=1,\dots,M, 
\edis 
without the use of the Jacob's ladders. 
\end{mydef13} 

\begin{mydefII1}
Since there is an infinite set of different choices of the sequence 
\bdis 
\{f_m(\sigma_0^m+it)\}_{m=1}^M
\edis 
from the set $\mfrak{D}$, we have also an infinite set of $\mfrak{D}$-equivalents of the Fermat-Wiles theorem. 
\end{mydefII1} 

\begin{remark}  
Of course, we can use, instead of the products (\ref{5.5}), also some linear combination of the integrals under consideration, and obtain new sets of the equivalents. 
\end{remark} 

\begin{remark}
If we define the $M$-segment 
\be \label{5.10} 
\mcal{F}_M=\prod_{m=1}^M\left[0,\FR\frac{\tau}{\left\{\prod_{m=1}^MF_m(\sigma_0^m)\right\}^{1/M}}\right], 
\ee  
then we can the $\mfrak{D}$-condition (\ref{5.9}) in the form 
\be \label{5.11} 
\lim_{\tau\to\infty}\frac{1}{\tau}\left\{
\underset{\mcal{F}_m}{\int\int\cdots\int}\prod_{m=1}^M|f_m(\sigma_0^m+it_m)|^2{\rm d}t_m
\right\}^{1/M}\not=1. 
\ee 
\end{remark} 

\section{On asymptotic connections between the set of Dirichlet's series, certain segments of Dirichlet's sum $D(x)$, Riemann's zeta-function and Jacob's ladders} 

\subsection{} 

Finally, we notice the following. The quotient of the formulas (\ref{3.10}) and (\ref{5.8}), factors $\frac{1}{\tau}$ are cancelled, gives the new asymptotic formula: 
\be \label{6.1} 
\begin{split}
& \left\{ 
\int_{\frac{x\tau}{\sqrt[3]{(1-c)\bar{c}(l)\zeta(2\sigma)}}}^{[\frac{x\tau}{\sqrt[3]{(1-c)\bar{c}(l)\zeta(2\sigma)}}]^1}\left|\zf\right|^2{\rm d}t\times \right. \\ 
& \left. \int_{1}^{\frac{x\tau}{\sqrt[3]{(1-c)\bar{c}(l)\zeta(2\sigma)}}}|\zeta(\sigma+it)|^2{\rm d}t\times \right. \\ 
& \left. \int_{1}^{\frac{x\tau}{\sqrt[3]{(1-c)\bar{c}(l)\zeta(2\sigma)}}}|S_1(t)|^{2l}{\rm d}t
\right\}^{1/3}\sim \\ 
& \prod_{m=1}^M\left\{
\int_0^{\frac{x}{\left\{\prod_{m=1}^MF_m(\sigma_0^m)\right\}^{1/M}}\tau}|f_m(\sigma_0^m+it)|^2{\rm d}t
\right\}^{1/M},\ \tau\to\infty 
\end{split}
\ee 
for every fixed 
\be \label{6.2} 
\begin{split}
& x>0,\ l,M\in\mbb{N},\ \sigma=\frac 12+\epsilon, \\ 
& f_m(\sigma_0^m+it)\in\mfrak{D},\ \sigma_0^m(f_m),\ m=1,\dots,M. 
\end{split}
\ee  

\subsection{} 

The quotient of the formula (\ref{3.10}), and, for example, the formula\footnote{See \cite{9}, (1.3) and (5.7).} 
\be \label{6.3} 
\lim_{\tau\to\infty}\frac{1}{\tau}\left\{
\sum_{\frac{x}{1-c}\tau<n\leq [\frac{x}{1-c}\tau]^1}d(n)
\right\}=x,\ x>0
\ee  
for some segments of the Dirichlet's sum $D(x)$, gives us the continuation of (\ref{6.1}): 
\be \label{6.4} 
\sim \sum_{\frac{x}{1-c}\tau<n\leq [\frac{x}{1-c}\tau]^1}d(n), 
\ee 
and so on. 

\begin{remark}
The formula (\ref{6.1}) together with its continuation (\ref{6.4}), expresses the asymptotic connection between the set of Dirichlet's series, some segments of the Dirichlet's sum $D(x)$, Riemann's zeta-function, and the Jacob's ladders. Of course, in our theory, there is an infinite set of possibilities for the next continuation of the initial chain (\ref{6.1}), (\ref{6.4}). 
\end{remark}

I would like to thank Michal Demetrian for his moral support of my study of Jacob's ladders.

\end{document}